\documentclass[12pt]{article}
\usepackage{amsmath, amsthm, amssymb}
\usepackage{graphicx}

\newcommand{\skein}[1]{\ensuremath{ {\mathcal S} ( {#1} )  }}
\newcommand{\solidtorus}{ \ensuremath{ S^1 \times D^2}}
\newcommand{\sumeven}[2]{ \sum_{\stackrel{k = #1}{ \text{even}}}^{#2} }

\newcommand{\IM}{\ensuremath{ I_A (M) }}

\newcommand{\firstroot}{ e ^{\pm 2 \pi i / 4r}}
\newcommand{\bracket}[2]{ \langle #1 \rangle_{#2}}
\newcommand{\symone}{ A^2 - A^{-2}}

\newtheorem{thm}{Theorem}
\newtheorem{cor}{Corollary}
\newtheorem{dfn}{Definition}
\newtheorem{lem}{Lemma}

\newcommand{\Proof}{\noindent {\bf Proof: \quad} }

\title{Quantum invariants can provide sharp Heegaard genus bounds}
\author{Helen Wong}
\date{}

\begin{document}

\maketitle
\author

\begin{abstract}
Using Seifert fibered three-manifold examples of Boileau and Zieschang, we demonstrate that
the Reshetikhin-Turaev quantum invariants may be used to provide a sharp lower
bound on the Heegaard genus which is strictly larger than the rank of the
fundamental group.
\end{abstract}

\section{Introduction}

For a closed orientable, connected three-manifold $M$, the Heegaard genus $g(M)$ is defined to be the smallest integer so that $M$ has a Heegaard splitting of that genus. 
Classically studied, the Heegaard genus is notoriously difficult to compute.  
However, in \cite{Garou} Garoufalidis showed that the  Reshetikhin-Turaev quantum invariants provide lower bounds on $g(M)$.
Though the quantum invariants may be computationally expensive, they allow an  algorithmic and combinatorial approach to studying the Heegaard genus. 

In this paper, we investigate one version of the Reshetikhin-Turaev quantum invariant which corresponds to the gauge group $SO(3)$.  
We show that under appropriate assumptions and choices, the calculations for the $SO(3)$ quantum invariants may be simplified.
Further, we give examples where the bounds to $g(M)$ are sharp. 

For a group $G$, let its rank $r(G)$ be the minimal number of elements required to generate $G$.  
The rank of the fundamental group of a three-manifold $r( \pi_1 M)$ is also a lower bound on $ g(M)$.  
By studying examples which first appeared in the work of Boileau and Zieschang (\cite{BZ}), we show that the $SO(3)$  quantum invariants may be used to provide a lower
bound on $g(M)$ which is strictly larger than $r( \pi_1 M)$.

The author would like to thank her thesis advisor Professor Andrew Casson for his patience and input while completing this research.

\section{Garoufalidis' lower bounds on Heegaard genus}

We begin by a brief review of the definition of the SO(3)-quantum invariant for $M$.  
In this paper, we follow the exposition of Lickorish in \cite{L2r} and of Turaev in \cite{T}, although the notation is changed slightly.

\begin{dfn}
Let $M$ be an oriented three-manifold and let $A$ be a complex number.
The skein space $\skein{M}$ is the complex vector space generated by all possible framed links and arcs in $M$, subject to the following Kauffman bracket relations:
\begin{itemize}
\item[i]  \quad Isotopy of framed links
\item[ii] \quad  \begin{minipage}{.3in}\includegraphics[width=.3in, , height=.45in]{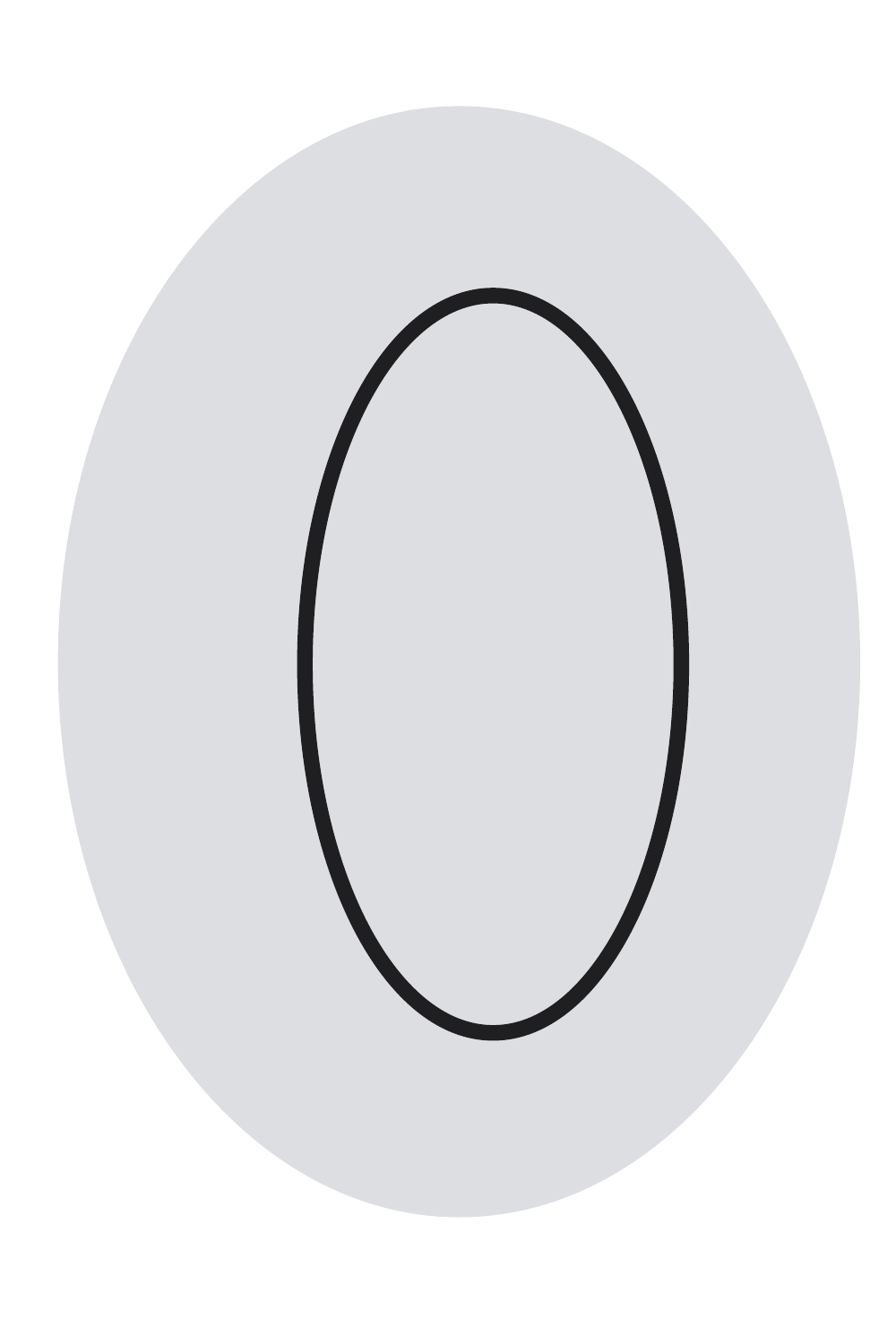} \end{minipage} = 
$ (A^2 - A^{-2}) $\begin{minipage}{.3in}\includegraphics[width=.3in, , height=.45in]{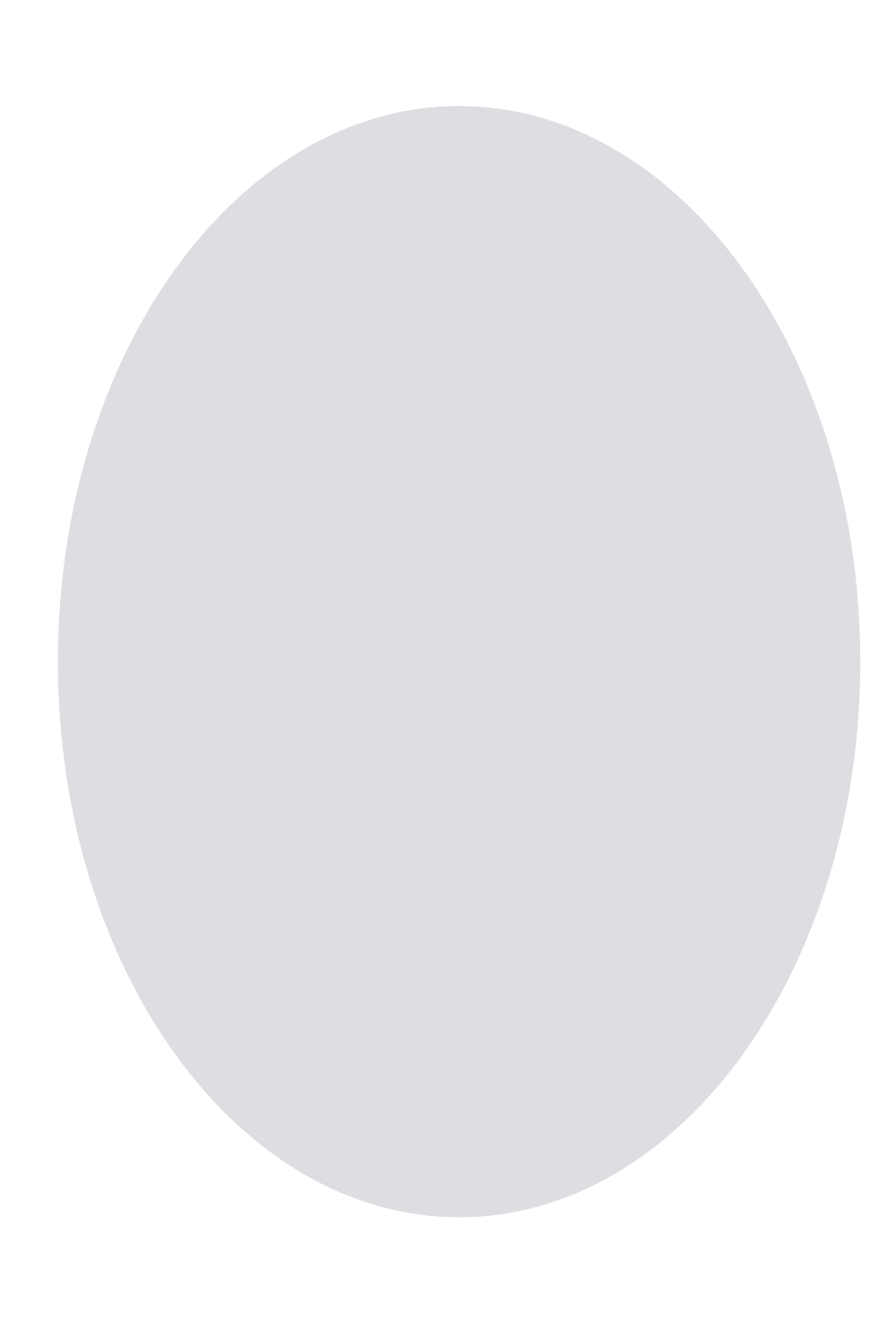} \end{minipage}
\item[iii] \quad \begin{minipage}{.3in}\includegraphics[width=.3in, , height=.4in]{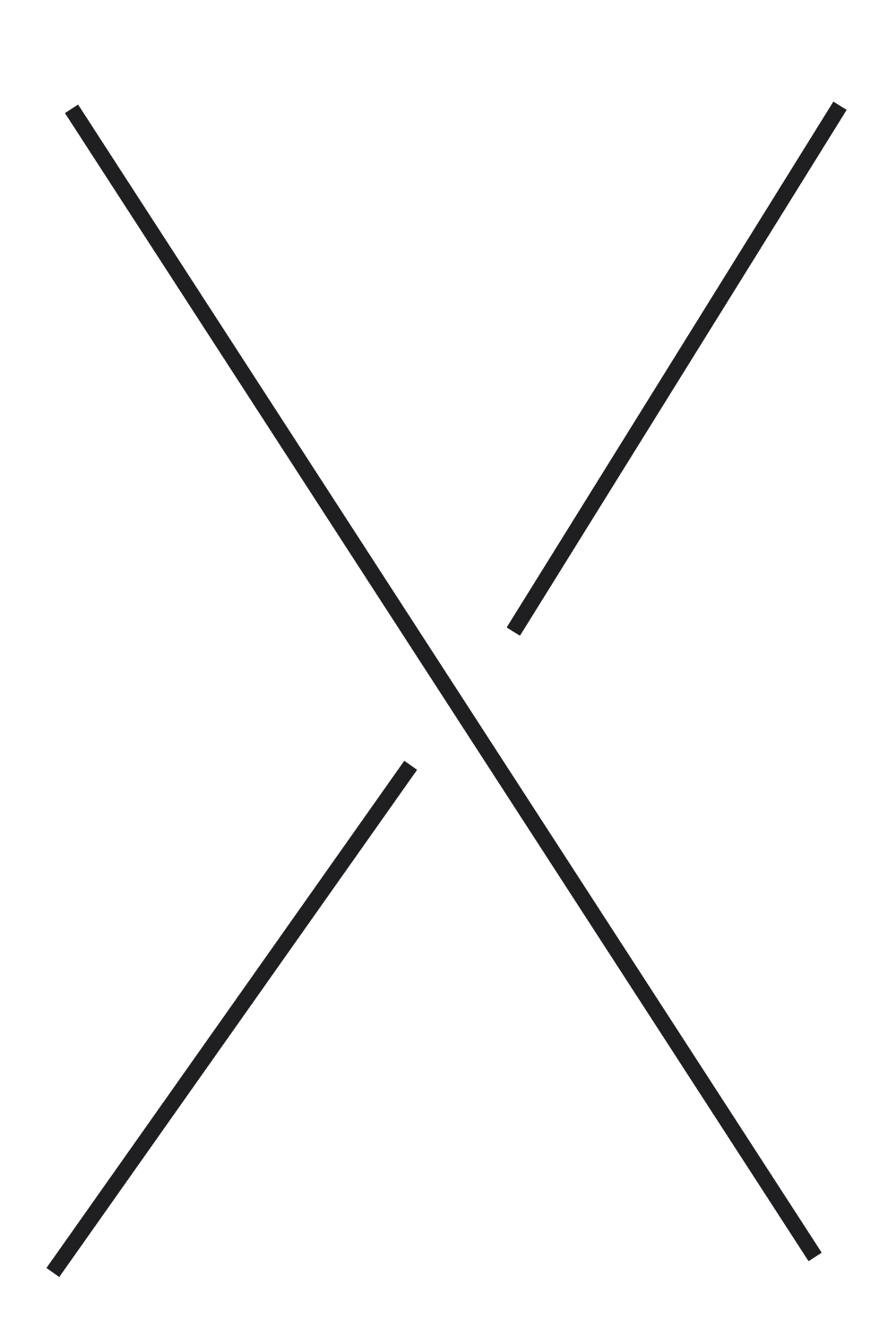} \end{minipage} 
$ = A $
\begin{minipage}{.3in}\includegraphics[width=.3in, , height=.4in]{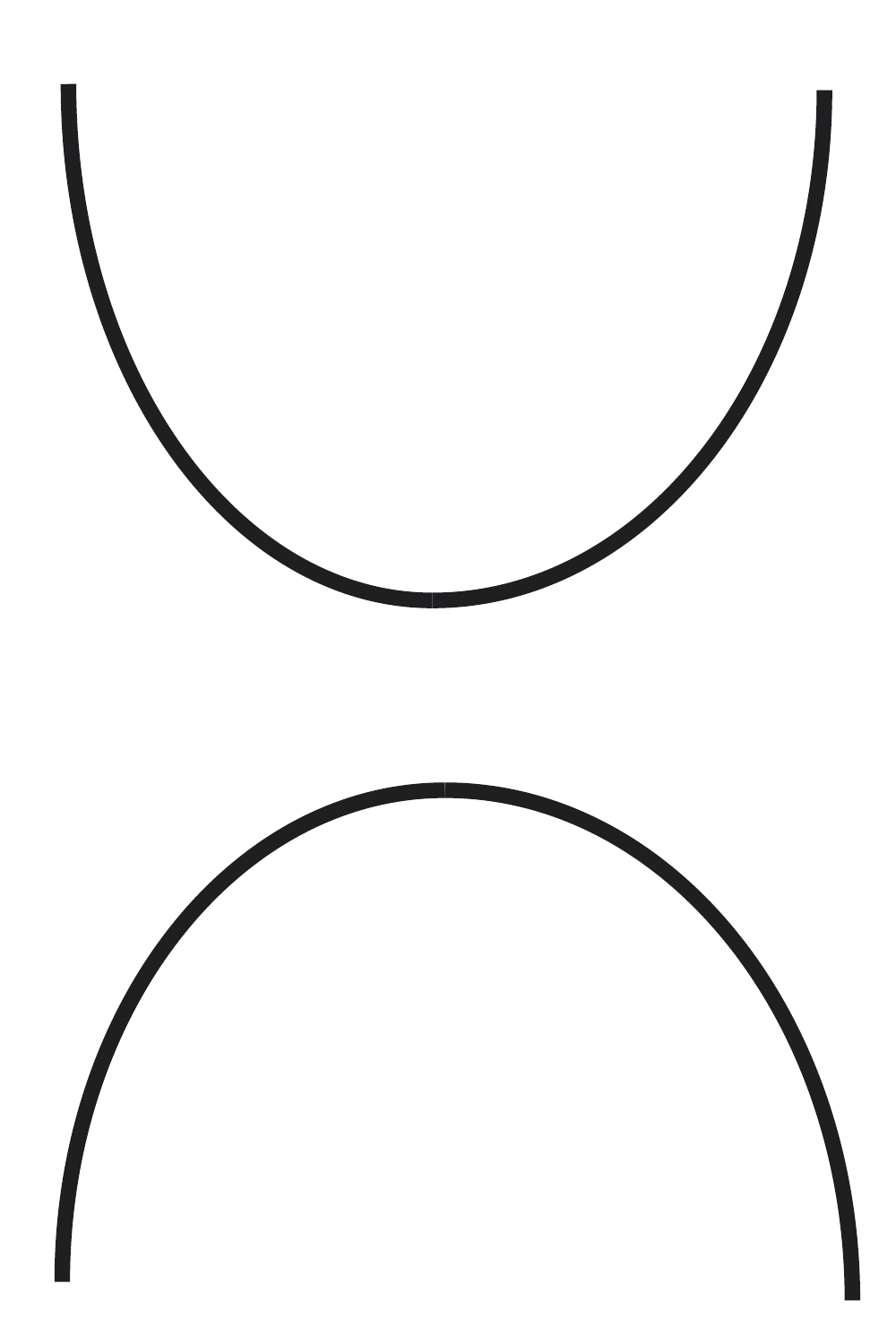} \end{minipage} 
$ + A^{-1} $
\begin{minipage}{.3in}\includegraphics[width=.3in, , height=.4in]{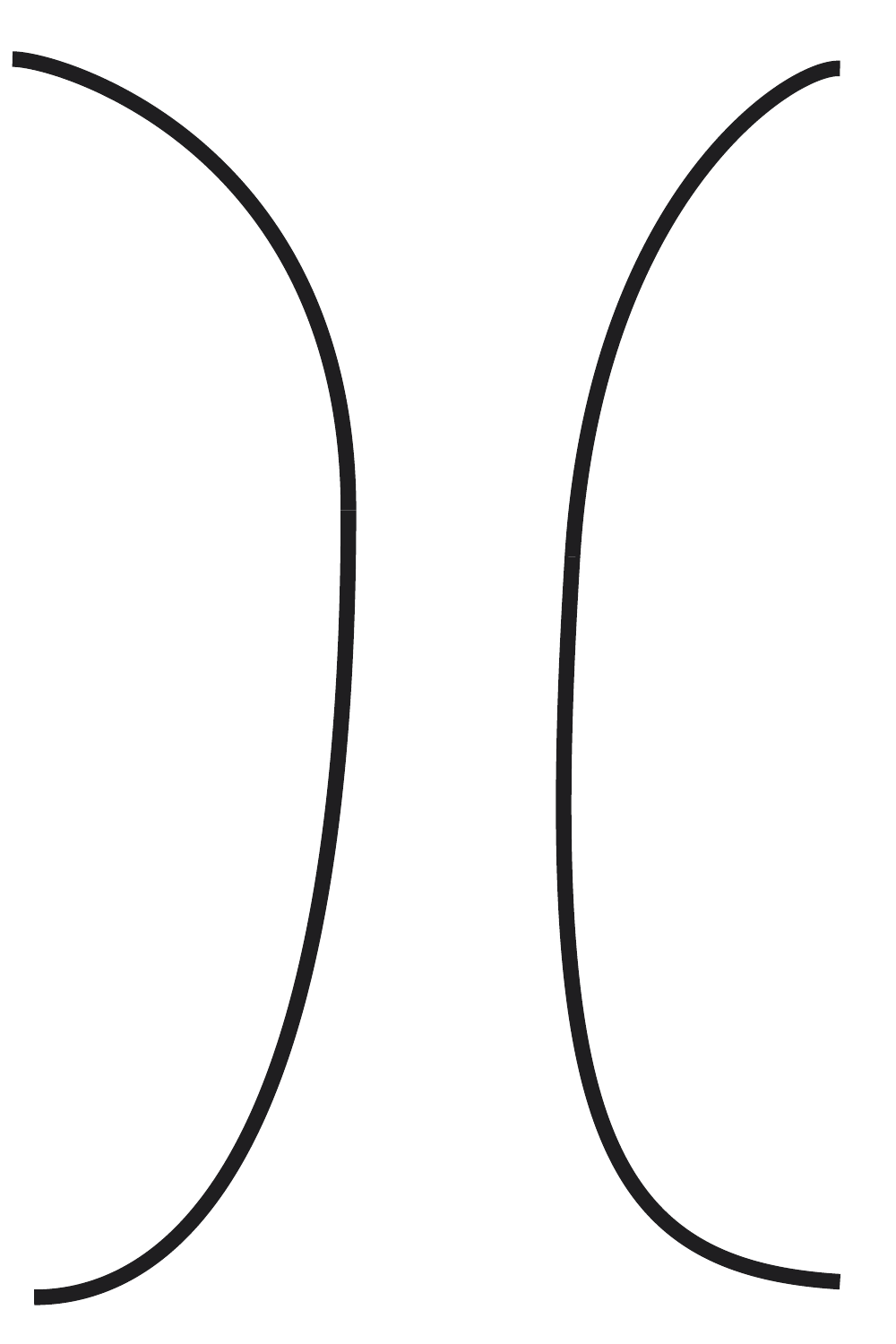} \end{minipage} .
\end{itemize}
\end{dfn}

The pictured diagrams are meant to describe a small ball near the framed link, outside of which the framed link is kept unchanged. 
For example,  it is easy to see that  $\skein{S^3} \cong \mathbb{C}$.  Also, $\skein{ \solidtorus }$ can be generated by the infinite set 
 $\mathcal B =  \Big \{   
 \begin{minipage}{.3in} \includegraphics[width=.3in, , height=.5in]{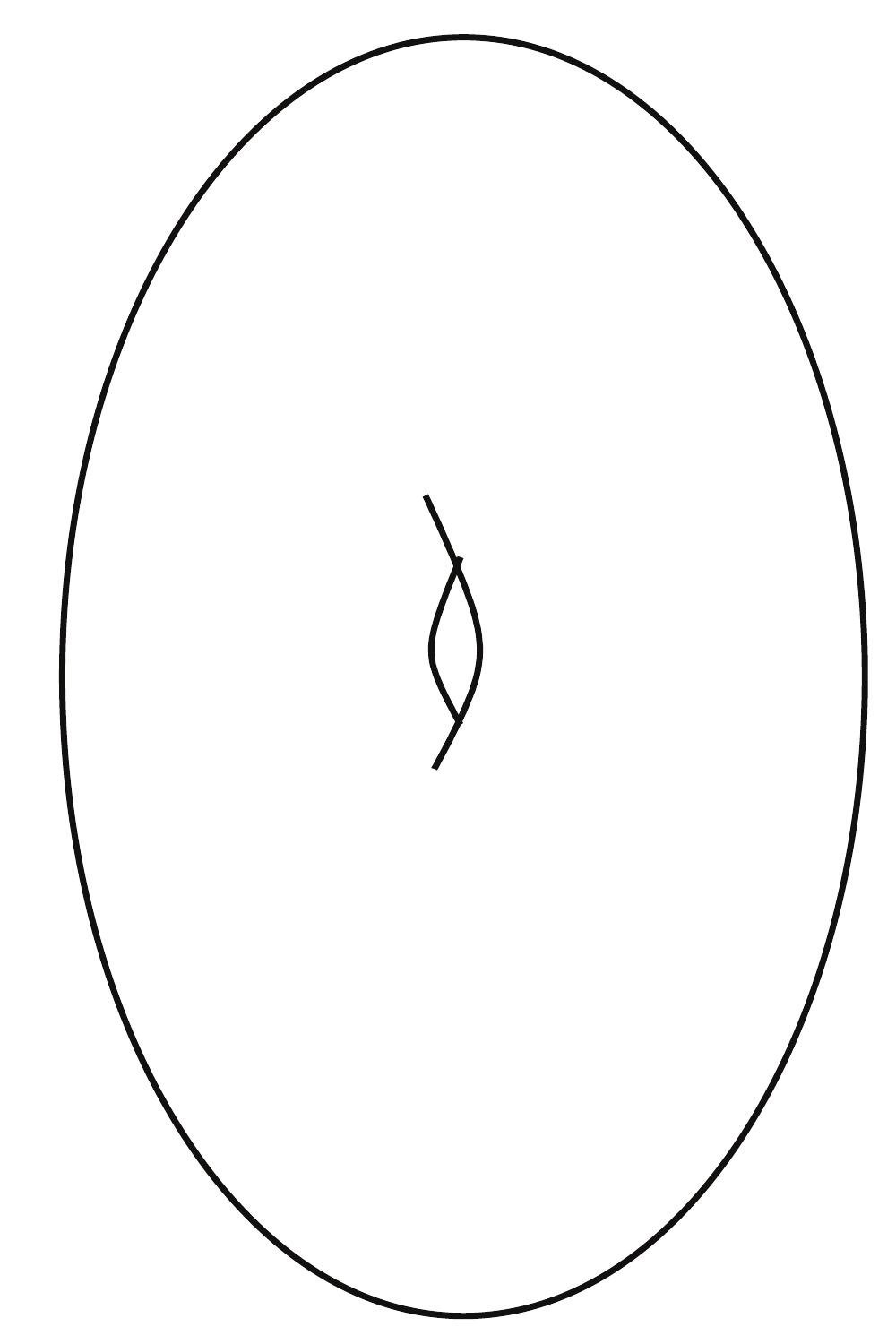} \end{minipage}, 
 \begin{minipage}{.3in} \includegraphics[width=.3in, , height=.5in]{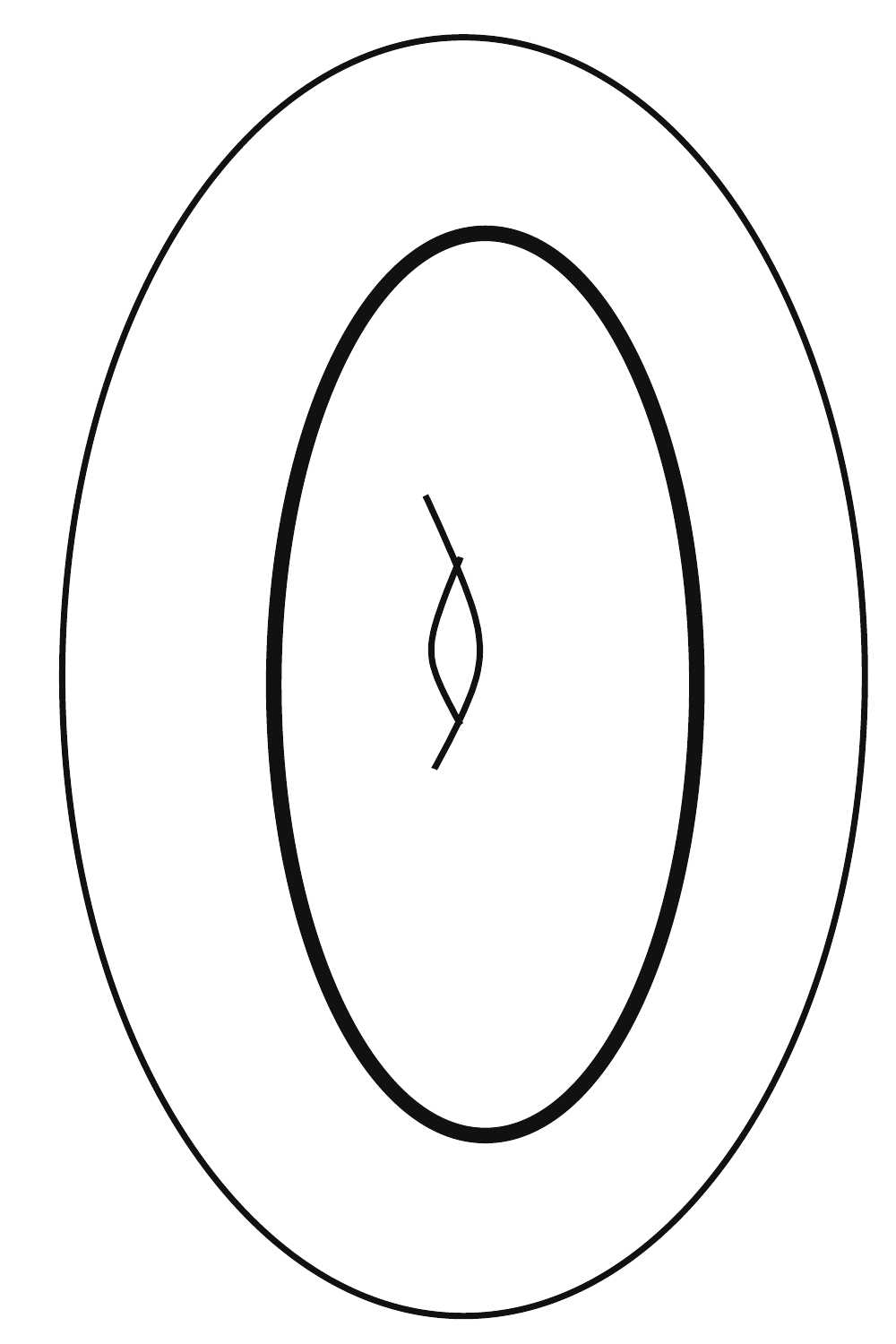} \end{minipage},
 \begin{minipage}{.3in} \includegraphics[width=.3in, , height=.5in]{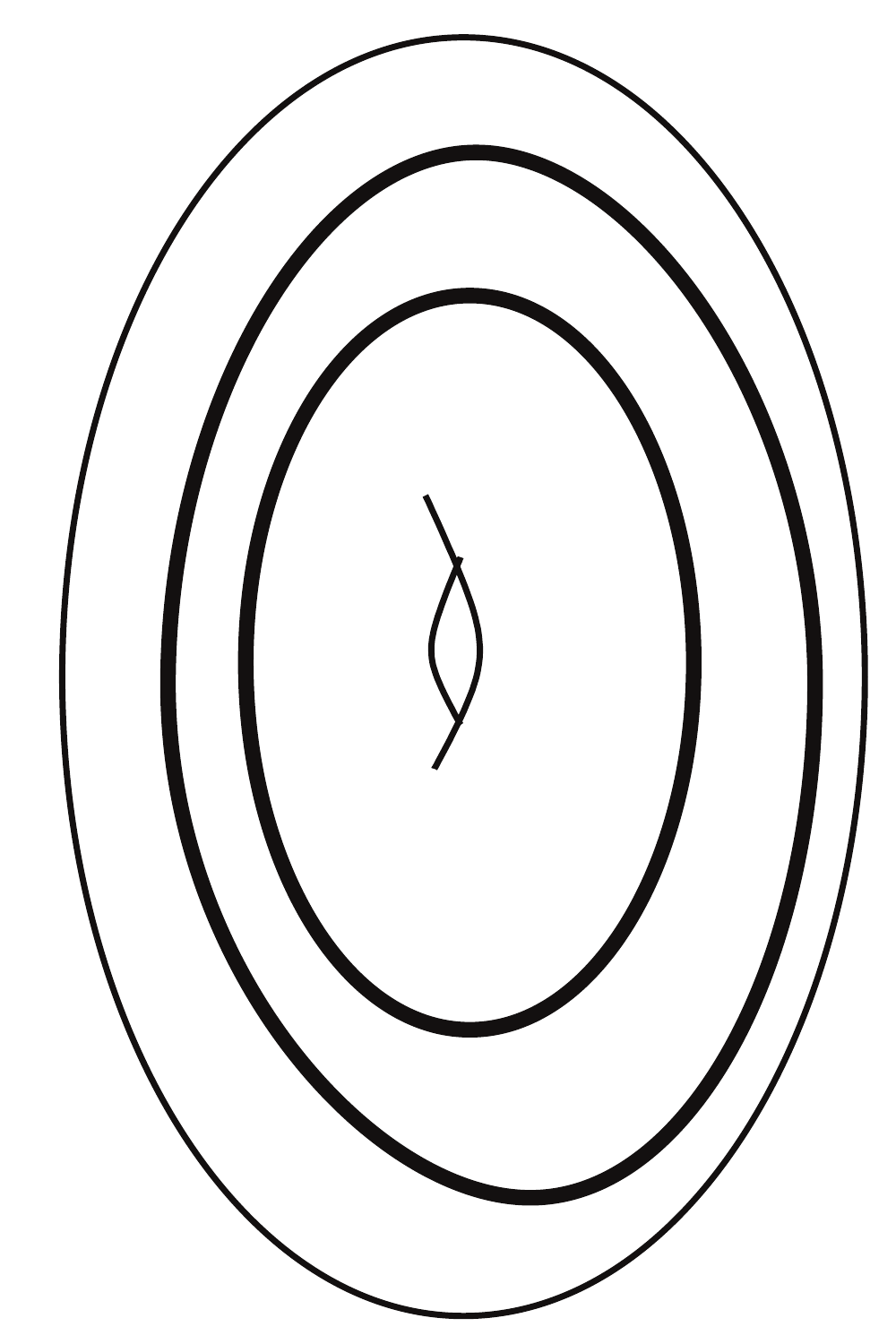} \end{minipage}, 
 \begin{minipage}{.3in} \includegraphics[width=.3in, , height=.5in]{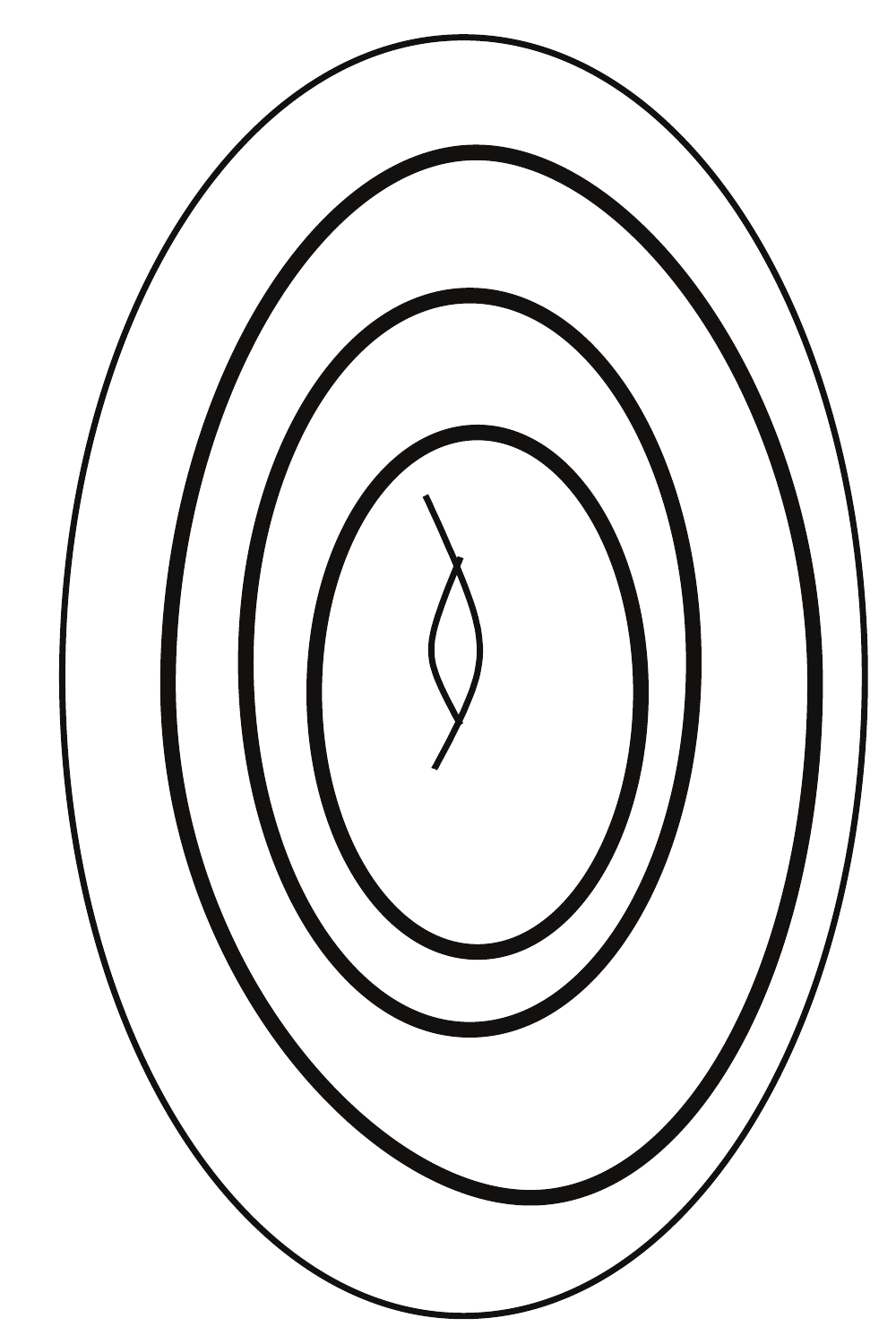} \end{minipage}, 
 \ldots
 \Big \}  $.

For a manifold $M$, choose a framed link $L$ in $S^3$ so that surgery along $L$ produces $M$.  
Work of Dehn and Lickorish \cite{LThm} guarantees such a link exists, and Kirby in \cite{KMoves} further shows that any such link is unique up to link isotopy and the two Kirby moves.   
As constructed by Lickorish \cite{L2r} and Turaev \cite{T}, the SO(3) quantum invariant for $M$ is obtained when each component of $L$ is replaced by the element $\Omega \in \skein{ \solidtorus }$. This element is carefully chosen so that, under a suitable normalization,  the resulting skein will be invariant under Kirby moves on $L$.  We describe $\Omega$  subsequently.  

To define the $SO(3)$ quantum invariants, we require $A$ to be a $2r$th or $4r$th primitive root of unity for some odd integer $r \geq 3$.  Then $\skein{ \solidtorus }$  becomes an $r$-dimensional vector space, and  there exists an alternative basis for $\skein{ \solidtorus }$ coming from Jones-Wenzl idempotents for the  group $SU(2)$.  We will denote the $k$-th basis element  by a open box labeled by $k$, drawn as follows:
$ S_k =  \begin{minipage}{.4in} \includegraphics[width=.4in, , height=.45in]{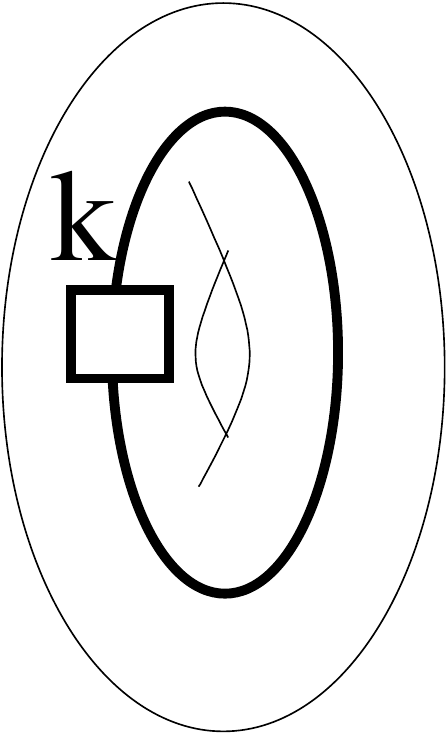} \end{minipage}$, with $k= 0, 1, \ldots, r-1
$.  
The $SO(3)$ quantum theory derives from consideration of the subspace generated by the evenly labelled Jones-Wenzl idempotents.

Let $\Omega $ be defined as follows, using only those Jones-Wenzl idempotents with even labels:
\[  \begin{minipage}{.4in} \includegraphics[width=.4in, , height=.45in]{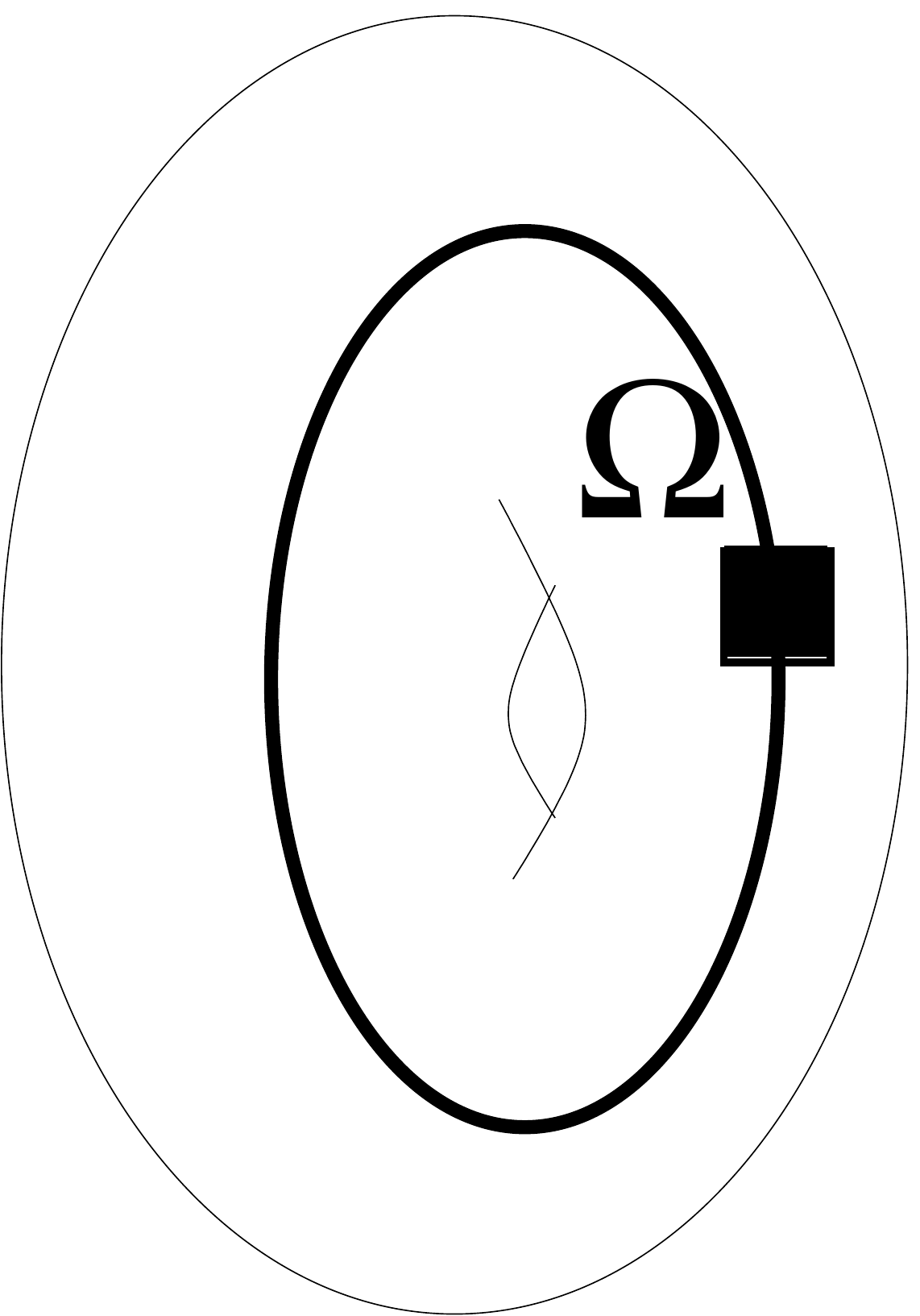} \end{minipage} =  \mu \sumeven{0}{r-3} \Delta_k   \; \begin{minipage}{.4in}\includegraphics[width=.4in, , height=.45in]{Skpicture.pdf}  \end{minipage} \;, \]
where  $ \mu ^2 = \frac{(A^2- A^{-2})^2}{-r} $ and $\Delta_k =  \frac{(-1)^k  \; ( A^{2(k+1)} -A^{-2(k+1)})}{ A^2 - A^{-2}}$.   
Let $\langle \Omega, \ldots, \Omega \rangle_L$ denote the skein in $S^3$ obtained by replacing each component of the link $L\subset S^3 $ by the skein element $\Omega$.  
Since $\skein{S^3} \cong \mathbb{C}$, we can reduce $\langle \Omega, \ldots, \Omega \rangle_L$ to a complex number dependent on $A$.

\begin{thm} (\cite{KM}, \cite{Blanchet}, \cite{L2r}) \label{Invariant}
Let $r \geq 3$ be an integer, and  $A$ be a $2r$th or $4r$th primitive root of unity.  
Let $M$ be the closed three-manifold which results from surgery along a framed link $L \subseteq S^3$, and let $\sigma(L)$ denote the signature of the linking matrix for $L$. \, Also let $U_-$ denote the unlink with framing $-1$ in $S^3$.  Then
\[ \IM = \langle \Omega, \ldots, \Omega \rangle_L   \left( \langle \Omega \rangle _{U_-} \right)^{\sigma(L)} \]
is an invariant of the three-manifold $M$.
\end{thm}

Theorem \ref{Invariant} appeared in various papers, notably ones by Kirby and Melvin \cite{KM} and Blanchet \cite{Blanchet}.  The version we state here is presented by Lickorish in \cite{L2r} and is referred to there as the ``invariant with zero spin structure".  In Turaev's book \cite{T}, it is the ``SO(3) quantum invariant".  

The SO(3) quantum invariant enjoys many properties; for instance, it behaves well under reversal of orientation and under connect sum.  That is, \[ I_A(\overline M) = \overline{\IM} \]
\[ I_A(M_1 \# M_2) = I_A(M_1) \# I_A(M_2) \] for three-manifolds $M, M_1, M_2$.

Further, for particular choices of $A$, it can be shown that $I_A(M)$ is related to the Heegaard genus, $g(M)$.  The proof may be found in \cite{Garou} and also \cite{T}.

\begin{thm} [\cite{Garou},  \cite{T}] 
Let $r\geq 3$ and $A= \firstroot $ or $ i \firstroot$.
Then 
$ | \IM| \leq  \mu^{-g(M)} $.
\end{thm}

Recall that $\mu^2 = \frac{ (A^2 - A^{-2})^2 }{-r}$.  When $A = \firstroot$  or $ i \firstroot$,  note that then $0 < \mu = 
\frac{2}{\sqrt{r}} \sin ( \frac{\pi}{r}) < 1$.
Define \[ q_A(M) = \log ( | \IM |) / \log(\mu) \]
so that $q_A(M) \leq g(M)$.  In other words, the SO(3) quantum invariant provides a lower bound on Heegaard genus.

\section{Changing the framing number by $r$}

We present some methods for simplifying the computation of the $SO(3)$ quantum invariants in some special cases.  In particular, it is possible for two non-homeomorphic manifolds to have $SO(3)$ quantum invariants with the same value.

\begin{thm} \label{Kr}
Let $r \geq 5$ be prime and $A$ be a $2r$th or $4r$th primitive root of unity.
Let $L$ be a link in $S^3$ with two distinct framings, $f$ and $f'$.  
Surgery along $L$ using the framings $f$ and $f'$ will result in two manifolds, which we call $M$ and $M'$ respectively.
Suppose that the framings $f$ and $f'$ are congruent modulo $r$ on each component of $L$.  Then \[ | I_A(M) | = | I_A(M') | \]
\end{thm}

\Proof
The proof follows from two lemmas about $\Omega$.  
\begin{lem}  The skein element $\Omega \in \skein{\solidtorus}$ does not change when any multiple of $r$ twists are added to it.
\[ \begin{minipage}{1in} \includegraphics[width=1in, , height=.75in]{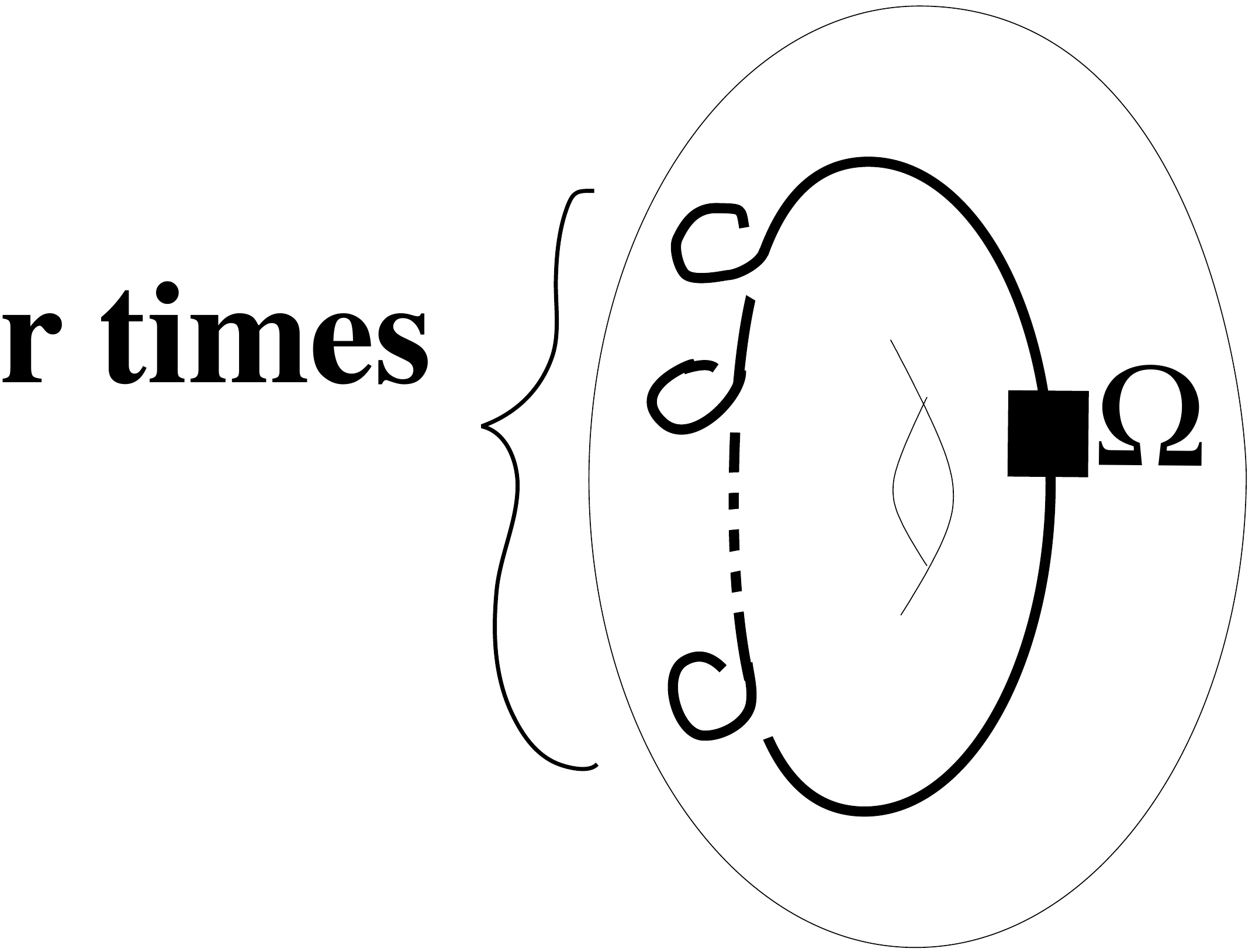} \end{minipage}
\; = \;
 \begin{minipage}{.6in} \includegraphics[width=.6in, , height=.75in]{Cw.pdf} \end{minipage}\]
\end{lem}
\Proof  This is an application of the following well-known fact about adding a twist to a Jones-Wenzl idempotent:
\[ \begin{minipage}{.45in} \includegraphics[width=.4in, , height=.5in]{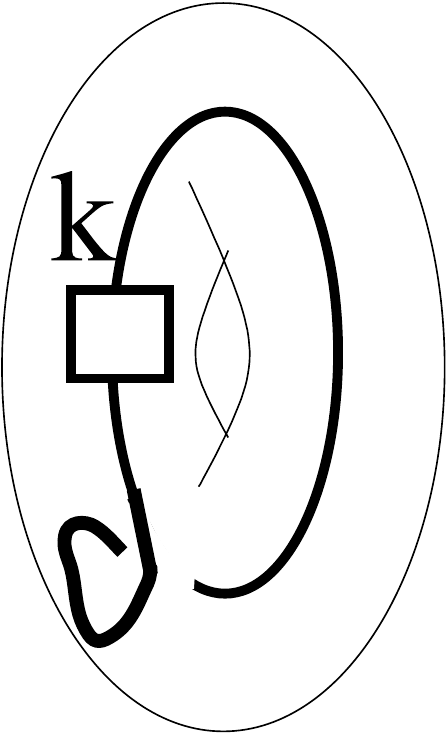} \end{minipage} 
 =  (-1)^k A^{k^2 + 2k}  \;
\begin{minipage}{.4in} \includegraphics[width=.4in, , height=.5in]{Skpicture.pdf} \end{minipage} \]
By definition $\Omega$ involves only labels $k$ which are even. 
Also
 $A^{4r}=1$ by assumption. 
 Thus we have that $\left((-1)^k A^{k^2+2k} \right)^r  = 1$, and the result follows.  
\\

The next lemma can be found in Lickorish, \cite{L2r}, or an article by Blanchet, \cite{Blanchet}.  The proof relies on calculation and standard facts about Gauss sums.  For details, see \cite{Methesis}.

\begin{lem} \label{omuminus}
\[ \bracket{\Omega}{U_-} = \mu \dfrac{- A^3 \left(  \sum_{k=0}^{r-1} A^{-4k^2} - (1/2) \sum_{k=0}^{4r-1} A^{-k^2} \right)}{ (\symone) }. \]
In particular, $|\bracket{\Omega}{U_-}  |= 1$.\\
\end{lem}
To finish the proof of Theorem \ref{Kr},  recall that a difference in framing numbers can be accounted for by introducing a corresponding number of twists into a diagram.
In particular, in the blackboard framing, the framing $f'$ can be obtained from $f$ by inserting $r$ twists, possibly more than once, to each link component.
Although this changes the signature of the framed link, Lemma \ref{omuminus} ensures that this will not affect the value of the quantum invariant upon taking aboslute value. $\Box$\\

Notice that Theorem \ref{Kr} is true only at the specified level $r$.  Examples of pairs of manifolds with all  values of the $SO(3)$ invariants identical for all choices of level $r$ can be found in  \cite{Lsame} and \cite{KB}.

We next recall a fact from number theory: any rational number $p/q \in \mathbb Q$ has a {\sl continued fraction decomposition}, where 
\[ \frac{p}{q} = x_0 - \cfrac{1}{ x_1 - \cfrac{1}{ \cdots - \frac{1}{x_n}}} .\]
We will denote this by $p/q = [x_0, x_1, \ldots, x_n]$.
We will say that two fractions $p/q$ and $p'/q'$ have {\sl entries in their continued fraction decompositions equal modulo $r$} if $p/q = [ x_0 ,  x_1, x_2, \ldots, x_m]$ and $p'/q' = [  x_0 + r a_0 ,  x_1+ r a_1, x_2 + r a_2 ,  \ldots, x_m + r a_m]$. 

Related by a series of Kirby moves, the following two surgery descriptions 
\[ \begin{minipage}{1.1 in} \includegraphics[width=1.1in, , height=.85in]{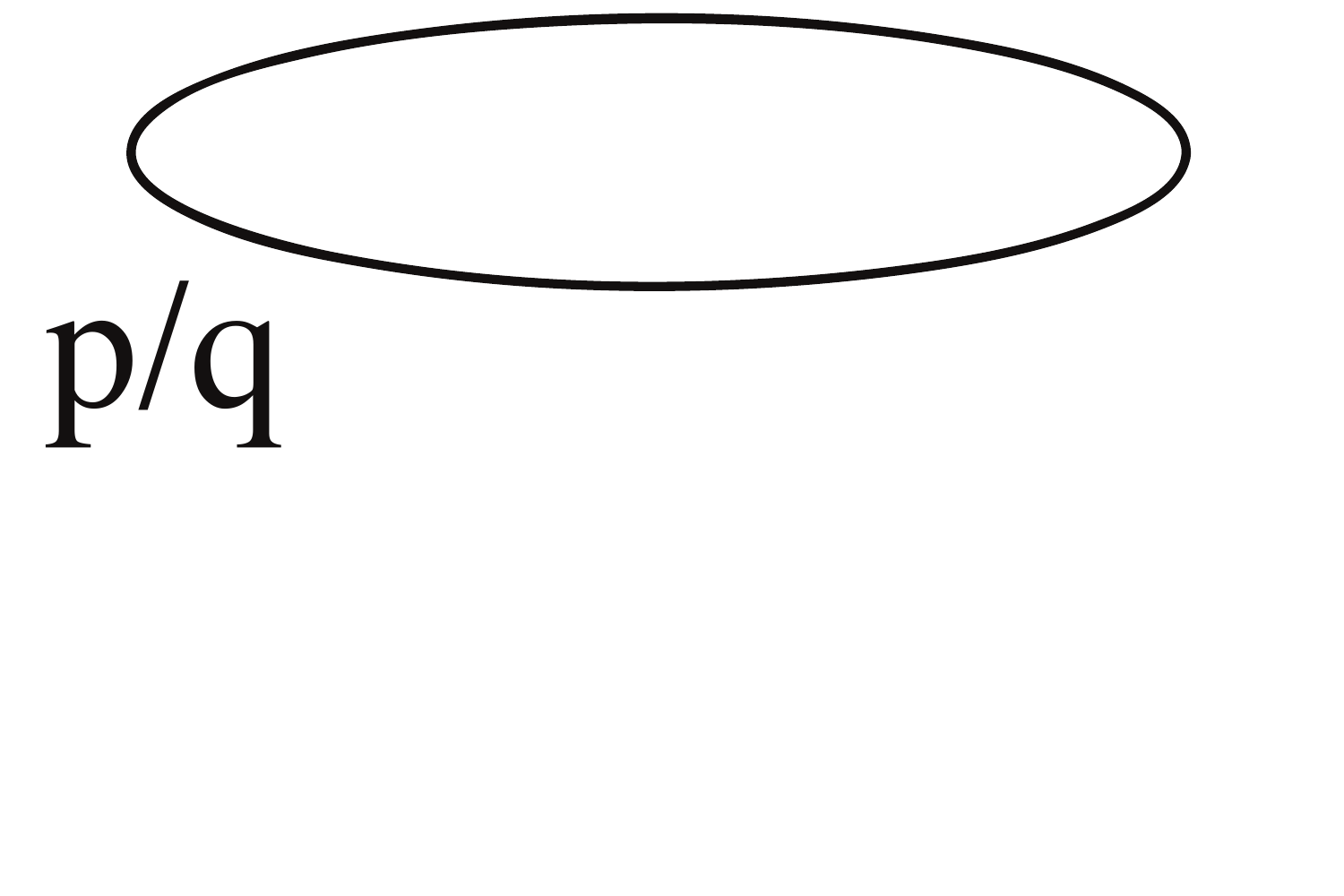} \end{minipage}
\; \text{  and  } \;
\begin{minipage}{1.1 in} \includegraphics[width=1.1in, , height=.85in]{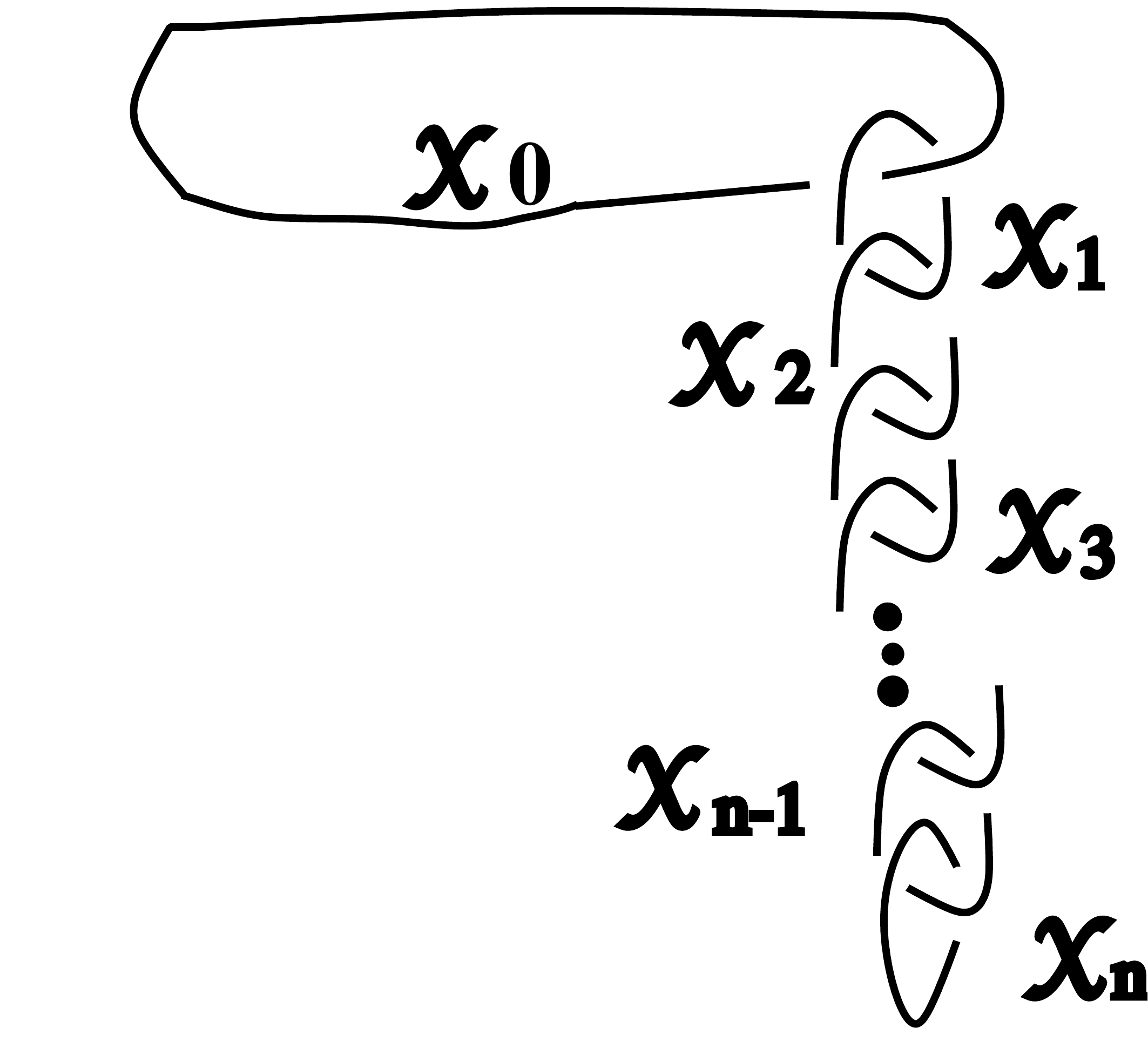} \end{minipage} \]
 yield the same manifold (\cite{Rolfsen}).
This allows us to convert a rational $(p/q)$-surgery along a knot into the language of surgery along a framed link, thus facilitating computation of the $SO(3)$ quantum invariants.

Let $L$ be a link with $\ell$ components in a three-manifold $M$.  Denote the manifold obtained by $(p_i / q_i)$-surgery along the $i$th component of $L$ ($1 \leq i \leq \ell$) by $M_{ \{p_i / q_i\} }$.  Theorem \ref{Kr} has the following corollary.
\begin{cor}  \label{cor}
Let $r \geq 5$ be prime and $A$ be a $2r$th or $4r$th primitive root of unity.
Let $L$ be a link in a three-manifold $M$.
If  $p_i/q_i $ and $p'_i/q'_i$ have entries in their continued fraction decomposition equal modulo $r$,
then \[ | I_A(M_{ \{p_i/q_i \}}) | = | I_A(M_{ \{p'/q'\}}) | \]
\end{cor}

\section{Boileau-Zieschang examples}

In this section, we focus on a particular set of three-manifolds.  
Let $M$ be the manifold which corresponds to the surgery presentation given by $L$ below:
\[ L = \begin{minipage}{1.5 in} \includegraphics[width=1.5in, , height=1in]{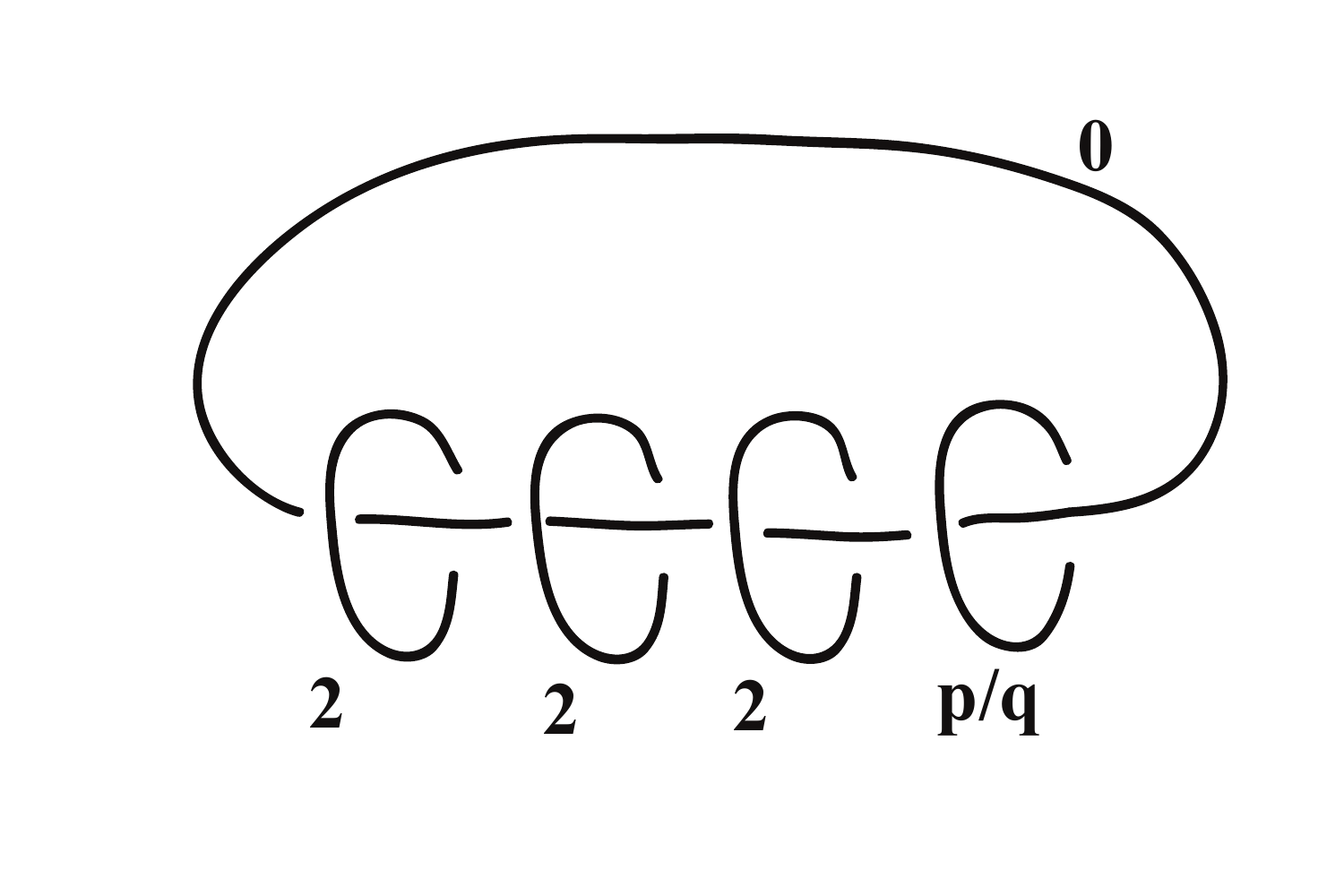}  \end{minipage} \]
where  $p/q$ haas every other entry iin its continued fraction decomposition divisible by $r$, i.e. $p/q = [r x_0,  x_1, r x_2, x_3, \ldots, x_{2n-1}, r x_{2n}]$ for  odd $r\geq 5$ and integers $x_i$. 

As defined here, $M$ is in fact a Seifert-fibered manifold and its Heegaard genus was first computed by Boileau and Zieschang in \cite{BZ}.  
Such $M$ are of especial interest because $2= r( \pi_1 M)$ and $g(M)=3$.  
For the specific values of $p/q$ chosen, we simplify the computation for the absolute value of the SO(3) quantum invariants by use of results in the previous section.  
We then compute that $2< q_A(M)$, and so
Garoufalidis' theorem implies that also  $2< g(M)$.
Thus viewed another way, we show that the quantum lower bound for Heegaard genus can be strictly larger than that provided by the rank of the fundamental group.

\begin{thm} \label{BZexample}
Let $r \geq 5$ be prime and $A$ be a $2r$th or $4r$th primitive root of unity .
Let $M$ be the manifold  corresponding to  surgery using the indicated coefficients along the link 
$L$ in $S^3$ described at the beginning of this section.
Suppose that $p/q = [r x_0,  x_1, r x_2, x_3, \ldots, x_{2n-1}, r x_{2n}]$ for some integers $x_{-1}, \ldots, x_{2n}$.
Then 
\[ | I_A(M) | = | I_A(  {\mathbb R}P^3 \# {\mathbb R}P^3 \# {\mathbb R}P^3 )|. \]
In particular,  $\mu^{-2}  < | I_A(M) | $ when $A= \firstroot$ or $i \firstroot$.
\end{thm}

\Proof
Let $U_{++}$ denote the the unlink with framing number $2$ in $S^3$.  (In the blackboard framing, this would be drawn as an unlink with two positive twists.)  Surgery along $U_{++}$ gives ${\mathbb R}P^3$, and by Lemma \ref{omuminus}, upon passing to absolute values we have $ | I_A( {\mathbb R}P^3) | = | \langle \Omega \rangle_{U_{++}} |$.

Because $p/q=[r x_0,  x_1, r x_2, x_3, \ldots, x_{2n-1}, r x_{2n}]$, the SO(3) invariant does not change absolute value when we replace the surgery label $p/q$ by $0= [0, x_1, 0, x_3, \ldots, x_{2n-1}, 0]$ according to Corollary \ref{cor}.  
Using Kirby moves, such a surgery picture is equivalent to the disjoint union of three copies of $U_{++}$.
Recall that if in a framed link, one component is an unknot with framing zero wich links only one other component geometrically once, then these two components may be deleted from the surgery picture without affecting any of the other remaining framed components.
Thus $| I_A(M) | = | I_A(  {\mathbb R}P^3 \# {\mathbb R}P^3 \# {\mathbb R}P^3 )|  =  | I_A( {\mathbb R}P^3) |^3 $.  

The next lemma again  follows from calculation using properties of the Jones-Wenzl idempotents.  We omit the proof and again refer to \cite{Methesis} for the explicit computation.

\begin{lem}
$ \langle \Omega \rangle_{U_{++}} = 
\mu \; \dfrac{ (A^{-4} - A^{2(r-1)})\sum^{r-1}_{k=0} A^{8k^2}}{ (\symone)^2} $.\\
\end{lem}
Assume that  $A = \cos( \frac{\pi}{ 2r}) + i \sin( \frac{\pi}{ 2r}) $ or 
$ i \cos( \frac{\pi}{ 2r}) -  \sin( \frac{\pi}{ 2r}) $.\\
Note that $\mu^2 = \frac{ (A^2 - A^{-2})^2 }{-r} = \frac{ 4 \sin^2(\pi/r)}{-r}$ in both cases. Also, $\sum^{r-1}_{k=0} A^{8k^2}  =  i^{ r(r-1)/2} \sqrt{r}$ from a standard result about Gauss sums, see for example \cite{BEW}.
It then follows that
\[ | I_A( {\mathbb R}P^3) | = | \langle \Omega \rangle_{U_{++}} |
 =  \dfrac{ \cos( \frac{\pi}{ 2r}) }{ \sin ( \frac{\pi}{ r}) } .\]
When $ r \geq 5$, a quick calculus argument confirms
\[ \mu^{-2} 
= \dfrac{ r }{4\sin^2(\frac{\pi}{r})} 
< \left( \dfrac{ \cos( \frac{\pi}{ 2r}) }{ \sin ( \frac{\pi}{ r}) } \right)^3 = | \langle \Omega \rangle_{U_{++}} |^3.\]
$\Box$

\begin{cor} \label{BZcor}
Let $M$ be the manifold  corresponding to  surgery using the indicated coefficients along the link $L$ pictured above in Theorem \ref{BZexample}.
Let $r \geq 5$ be odd, and let $ A = \firstroot $ or $i \firstroot$.
Suppose also that $p/q = [r x_0,  x_1, r x_2, x_3, \ldots, x_{2n-1}, r x_{2n}]$ for some integers $x_{-1}, \ldots, x_{2n}$. 
Then  \[2< q_A(M) \leq g(M).\]
\end{cor}

\Proof 
This follows directly from Theorem \ref{BZexample} and the definition of $q_A(M)$, the lower bound on Heegaard genus provided by Garoufalidis' Theorem.
$\Box$



\begin{thebibliography}{99}
\addcontentsline{toc}{chapter}{\numberline{}Bibliography}





\bibitem{BEW} Berndt, B.; Evans, R.; and Williams, K.,  Gauss and Jacobi sums. Canadian Mathematical Society Series of Monographs and Advanced Texts. A Wiley-Interscience Publication. John Wiley and Sons, Inc., New York, 1998.


\bibitem{Blanchet} Blanchet, C., {\sl Invariants on three-manifolds with spin structure}. Comment. Math. Helv. 1992 {\bf 67}, no. 3, 406-427.



\bibitem{BZ} Boileau, M. and Zieschang, H., {\sl Heegaard genus of closed orientable Seifert 3-manifolds}. Invent. Math. 1987 {\bf 76}, 455-468.






\bibitem{Garou} Garoufalidis, S., {\sl Applications of quantum invariants in low-dimensional topology}.  Topology, 1998  {\bf 37}, 219--224.


\bibitem{KB} Kania-Bartoszynska, J., {\sl Examples of different $3$-manifolds with the same invariants of Witten and Reshetikhin-Turaev}. 
Topology 1993 {\bf 32},  47--54. 

\bibitem{KMoves} Kirby, R., {\sl A calculus for framed links in $S^3$}. Invent. Math. 1978 {\bf 45}, 36-56.

\bibitem{KM} Kirby, R. and Melvin, P., {\sl The 3-manifold invariants of Witten and Reshetikhin-Turaev for sl(2, C)}.  Invent. Math. 1991, {\bf 105}, 473-545.




\bibitem{LSurvey} Lickorish, W. B. R.,  {\sl Quantum invariants of 3-manifolds}.  Handbook of geometric topology,  707--734, North-Holland, Amsterdam, 2002.

\bibitem{LThm} Lickorish, W.B. R., {\sl A representation of orientable combinatorial 3-manifolds}. Ann. Math. 1962 {\bf 76}, no.2, 531-540.

\bibitem{L2r} Lickorish, W. B. R., {\sl The skein method for three-manifold invariants}. J. Knot Theory Ramifications 1993 {\bf} 2, no. 2, 171-194.

\bibitem{Lsame} Lickorish, W. B. R., {\sl 
Distinct $3$-manifolds with all ${\rm SU}(2)\sb q$ invariants the same }. 
Proc. Amer. Math. Soc.  1993 {\bf 117}, 285--292. 










	\bibitem{RT} Reshetikhin, R. and Turaev, V., {Invariants of 3-manifolds via link polynomials and quantum groups}.  Invent. Math. 1991, {\bf 103}, 547-597.




\bibitem{Rolfsen}  Rolfsen, D., Knots and links. Mathematics Lecture Series, {\bf  7}. Publish or Perish, Berkeley, 1976.


 
\bibitem{T} Turaev, V. G.,   Quantum invariants of knots and 3-manifolds. de Gruyter Studies in Mathematics,  {\bf 18} . Walter de Gruyter \& Co., Berlin, 1994. 

\bibitem{Waldhausen} Waldhausen, F., {\sl Heegaard-Zerlegungen der 3-spare}. Topology, 1968 {\bf 7}, 196-203.

\bibitem{Methesis} Wong, H.,  The SO(3) quantum invariants:  Their density and topological applications.  Ph. D. thesis, Yale University, 2007.


\end{thebibliography}
\end{document}